\documentclass[12pt,report]{amsart}
\usepackage{graphicx}
\usepackage[all,cmtip]{xy}
\usepackage[margin=1in]{geometry}
\usepackage{amssymb}
\usepackage{amsthm}
\usepackage{hyperref}
\hypersetup{pdfstartview={FitH}}
\usepackage{multicol}
\usepackage{latexsym}
\usepackage{amsfonts}
\usepackage{amscd}
\usepackage[all]{xy}
\usepackage{slashed}
\usepackage{todonotes}
\let\oldtocsection=\tocsection

\let\oldtocsubsection=\tocsubsection

\let\oldtocsubsubsection=\tocsubsubsection

\renewcommand{\tocsection}[2]{\hspace{0em}\oldtocsection{#1}{#2}}
\renewcommand{\tocsubsection}[2]{\hspace{1em}\oldtocsubsection{#1}{#2}}
\renewcommand{\tocsubsubsection}[2]{\hspace{2em}\oldtocsubsubsection{#1}{#2}}

\newtheorem{proposition}{Proposition}
\newtheorem{theorem}[proposition]{Theorem}
\newtheorem{corollary}[proposition]{Corollary}

\newtheorem{thm}[proposition]{Theorem}

\newtheorem{define}[proposition]{Definition}
\newtheorem{rem}[proposition]{Remark}

\newcommand{\Z}{\mathbb{Z}}
\newcommand{\R}{\mathbb{R}}
\newcommand{\C}{\mathbb{C}}

\newcommand{\Q}{\mathbb{Q}}
\newcommand{\U}{\mathbb{U}}
\newcommand{\F}{\mathbb{F}}

\makeatletter
\newcommand{\xyC}[1]{%
\makeatletter
\xydef@\xymatrixcolsep@{#1}
\makeatother
} % end of \xyR
\makeatletter
\newcommand{\xyR}[1]{%
\makeatletter
\xydef@\xymatrixrowsep@{#1}
\makeatother
} % end of \xyR

\def\XXint#1#2#3{{\setbox0=\hbox{$#1{#2#3}{\int}$}
\vcenter{\hbox{$#2#3$}}\kern-.5\wd0}}

\title{Abelian Analytic Torsion and Symplectic Volume}

\author{B.D.K. McLellan}
\address{Northeastern University, Department of Mathematics, 360 Huntington Ave, Boston, MA, USA, 02115}
\email{b.mclellan@neu.edu}

\begin{document}

\maketitle

\begin{abstract}
This article studies the abelian analytic torsion on a closed, oriented, Sasakian three-manifold and identifies this quantity as a specific multiple of the natural unit symplectic volume form on the moduli space of flat abelian connections.  This identification computes the analytic torsion explicitly in terms of Seifert data.
\end{abstract}

\tableofcontents

%\begin{keyword}
%analytic torsion, contact torsion, Chern-Simons theory, Sasakian three-manifold, quantum field theory
%\end{keyword}

\section{Introduction}\label{nonabelchap}
\noindent
This article studies the abelian analytic torsion on Sasakian three-manifolds.  The analytic torsion is a topological invariant that was introduced by D.B. Ray and I.M. Singer \cite{rsi} as an analytic analogue of the combinatorially defined Reidemeister torsion \cite{reid2}.  It is a well known fact that these two torsions agree, as was independently shown by W. M\"{u}ller, \cite{mul}, and J. Cheeger, \cite{cheeg}, for unimodular representations.  More recently an elegant new proof of this equivalence has been given by M. Braverman \cite{brav} using the Witten laplacian \cite{wit}.\\
\\
Our main objective in this article is to compute the (square-root of the) analytic torsion explicitly as a natural symplectic volume form on the moduli space of flat abelian connections.  This identification is motivated by the work of C. Beasley and E. Witten \cite{bw} involving Chern-Simons theory on contact three-manifolds.  Recall that A.S. Schwarz \cite{schw} has shown that the abelian Chern-Simons partition function is proportional to the analytic torsion and our study is also natural in light of this fact.  Our main result, Theorem \ref{maintheorem}, shows that two mathematically a priori different definitions of the abelian Chern-Simons partitition function derived from \cite{bw} are rigorously equivalent.  Our main strategy is to use the the work of M. Rumin and N. Seshadri \cite{rs} which naturally connects the analytic torsion with contact structures on three-manifolds.  
\\
\\
Throughout, $X$ will denote a closed, orientable three-manifold, and $(X,\phi,\xi,\kappa,\operatorname{G})$ will denote $X$ equipped with a \emph{Sasakian} structure.  See \cite{bg}, \cite{dt} for standard background on Sasakian and contact geometry.  For convenience we recall that a \emph{Sasakian manifold} is a normal contact metric manifold, $(X,\phi,\xi,\kappa,\operatorname{G})$, where,
\begin{itemize}
\item  $\kappa\in\Omega^{1}(X)$ is a contact form, i.e. $\kappa\wedge d\kappa \neq 0$, $\xi\in\Gamma(TX)$ is the Reeb vector field,

\item  $\phi\in \operatorname{End}(TX)$, $\phi(Y)=:JY$ for $Y\in \Gamma(H)$, $\phi(\xi)=0$ where $J\in \operatorname{End}(H)$ is an almost complex structure on the contact distribution $H:=\operatorname{ker}\,\kappa\subset TX$, and,

\item  $\operatorname{G}=\kappa\otimes\kappa+d\kappa\circ (\mathbb{I}\otimes \phi)$.
\end{itemize}
\begin{define}
A \emph{Seifert manifold} is a closed orientable three-manifold that admits a locally free $\U(1)$-action.
\end{define}
\begin{rem}
See \cite{o1} for a general definition and classification of Seifert manifolds.
\end{rem}
\noindent
Let $\Sigma$ denote the base of a Seifert manifold when viewed as the total space of a $\U(1)$-bundle,
\begin{displaymath}
\xymatrix{\xyC{2pc}\xyR{1pc}\U(1) \ar@{^{(}->}[r] & X \ar[d]\\
                              & \Sigma}.
\end{displaymath}
It is well known that the topological isomorphism class of a Seifert manifold $X$ is determined by its Seifert invariants \cite{o1},
\begin{equation*}
[g,n; (\alpha_1,\beta_1),\ldots,(\alpha_M,\beta_M)],\,\,\, \operatorname{gcd}(\alpha_j,\beta_j)=1,
\end{equation*}
where $g$ is the genus of $\Sigma$.  Geometrically, the $\U(1)$ action on $X$ is rotations of the fibres over $\Sigma$ and the points in the $\U(1)$ fiber over each orbifold point $p_j$ on $\Sigma$ are fixed by the cyclic subgroup $\Z_{\alpha_j}$ of $\U(1)$.  The fundamental group $\pi_{1}(X)$ is generated by the following elements \cite{o1},
\begin{eqnarray*}
& a_p,b_p,\,\,\, p=1,\ldots,g,\\
& c_j, \,\,\, j=1,\ldots,M,\\
& h,
\end{eqnarray*}
which satisfy the relations,
\begin{eqnarray}\label{fundpres}
[a_p,h]=[b_p,h]=[c_j,h]&=&1,\\
c_j^{\alpha_j}h^{\beta_j}&=&1\nonumber,\\
\prod_{p=1}^{g}[a_p,b_p]\prod_{j=1}^{M}c_j&=&h^n,\nonumber
\end{eqnarray}
Geometrically, the generator $h$ is associated to the generic $\U(1)$ fiber over $\Sigma$, the generators $a_p, b_p$ come from the $2g$ non-contractible cycles on $\Sigma$, and the generators $c_j$ come from the small one cycles in $\Sigma$ around each of the orbifold points $p_j$.  
\begin{rem}
Since the analytic torsion is defined with respect to a choice of metric, we naturally work with Sasakian structures.  Recall that $X$ admits a Sasakian structure $(X,\phi,\xi,\kappa,\operatorname{G})$ $\iff$
\begin{itemize}
\item  \cite[Theorem 7.5.1, 7.5.2]{bg} $X$ admits a Seifert structure that is the total space of a non-trivial principal $\U(1)$ orbibundle over a Hodge orbifold surface, $\Sigma$.
\end{itemize}
For this article, Seifert structures on $X$ are induced by Sasakian structures.  
\end{rem}
\noindent
Let $\mathbb{T}$ denote a compact, connected abelian Lie group of real dimension $N$, $\frak{t}$ denote its Lie algebra and $\Lambda\subset \frak{t}$ the integral lattice.  Let $\operatorname{Tors}H^{2}(X,\Lambda)$ denote the torsion subgroup of $H^{2}(X,\Lambda)$.  For $P$ a principal $\mathbb{T}$-bundle over $X$, $\mathcal{A}_{P}$ is the affine space of connections on $P$ modeled on the vector space of $\mathbb{T}$-invariant horizontal one-forms on $P$, $(\Omega^{1}_{\operatorname{hor}}(P,\frak{t}))^{\mathbb{T}}\simeq\Omega^{1}(X,\frak{t})$.  The group of smooth gauge transformations is the group of $\mathbb{T}$ equivariant smooth maps $\mathcal{G}:=(\operatorname{Map}^{\infty}(P,\mathbb{T}))^{\mathbb{T}}\simeq \operatorname{Map}^{\infty}(X,\mathbb{T})$ and acts on $\mathcal{A}_{P}$ in the standard way.  That is, for $g\in\operatorname{Map}^{\infty}(P,\mathbb{T})$, and $A\in \mathcal{A}_{P}$, $A\cdot g:=A+g^{*}\vartheta$, where $\vartheta\in \Omega^{1}(\mathbb{T},\frak{t})$ denotes the Maurer-Cartan form on $\mathbb{T}$.  In order to define the Chern-Simons action, a negative definite symmetric bilinear form on $\frak{t}$ needs to be chosen.  Let $B\mathbb{T}$ denote the classifying space of principal $\mathbb{T}$-bundles.  Valid choices for such forms $\langle\cdot, \cdot\rangle\in\operatorname{Sym}^{2}_{\mathbb{T}}(\frak{t}^{*})$ are classified by elements of $H^{4}(B\mathbb{T},\mathbb{Z})$ \cite{dw1}, \cite{bm}.  Choosing a basis $e^{\alpha}$ for $H^{2}(B\mathbb{T},\mathbb{Z})$, an element in $H^{4}(B\mathbb{T},\mathbb{Z})$ may be written as $M_{\alpha\,\beta}\, e^{\alpha}\cup e^{\beta}$, where $M_{\alpha\, \beta}$ is an $N\times N$ integral symmetric matrix.  For the purposes of this article we choose $\langle\cdot, \cdot\rangle$ corresponding to $M_{\alpha\,\beta}=-2\mathbb{I}_{\alpha\,\beta}$, where $\mathbb{I}_{\alpha\,\beta}$ is the identity matrix.  Let $W$ be a compact oriented four-manifold such that $\partial W=X$, which always exists \cite{ps}.  Extend $P$ to a $\mathbb{T}$-bundle $Q$ over $W$, which is always possible in our case \cite{bm}.  Given a form $\alpha\in\Omega^{j}(P,\frak{t})$, let $\widetilde{\alpha}\in \Omega^{j}(Q,\frak{t})$ denote the corresponding extension to $Q$.  For a connection $A\in\Omega^{1}(P,\frak{t})$, denote the curvature form of the extension $\widetilde{A}\in\Omega^{1}(Q,\frak{t})$ by $F_{\widetilde{A}}\in\Omega^{2}(W,\frak{t})$.
\begin{define}\label{csactiondef}
The \emph{Chern-Simons action} of a $\mathbb{T}$-connection $A\in\mathcal{A}_{P}$ is defined by,
\begin{equation}\label{act}
\operatorname{CS}_{X,P}(A):=\frac{1}{4\pi}\int_{W}\langle F_{\widetilde{A}}\wedge F_{\widetilde{A}}\rangle\,\,\,\operatorname{mod}\,\, (2\pi\Z).
\end{equation}
\end{define}
We also define the following,
\begin{itemize}
\item  $m_X:=\frac{N}{2}(\operatorname{dim}H^{1}(X;\R)-2\operatorname{dim}H^{0}(X;\R))$,

\item  $A_{P}$ denotes a flat connection on a principal $\mathbb{T}$-bundle $P$ over $X$,

%\item   $[g,n; (\alpha_{1},\beta_{1}),\ldots,(\alpha_{M},\beta_{M})]$, for $\operatorname{gcd}(\alpha_{j},\beta_{j})=1$) are the Seifert invariants of a Sasakian manifold $(X,\phi,\xi,\kappa,\operatorname{G})$,

\item  $c_{1}(X)=n+\sum_{j=1}^{M}\frac{\beta_{j}}{\alpha_{j}}$ is the first orbifold Chern number of the Seifert manifold $X$,

\item  $s(\alpha,\beta):=\frac{1}{4\alpha}\sum_{j=1}^{\alpha-1}\operatorname{cot}\left(\frac{\pi j}{\alpha}\right)\operatorname{cot}\left(\frac{\pi j\beta}{\alpha}\right)\in\Q$ is the Rademacher-Dedekind sum,

\item  $\eta_{0}=N\left(\frac{c_{1}(X)}{6}-2\sum_{j=1}^{M}s(\alpha_{j},\beta_{j})\right)$ is the \emph{adiabatic eta-invariant} of the Sasakian manifold $(X,\phi,\xi,\kappa,\operatorname{G})$ \cite{nic},

\item  $\mathcal{M}_{X}\simeq \coprod_{[P]\in\operatorname{Tors}H^{2}(X,\Lambda)}\mathbb{T}^{2g}$ denotes the moduli space of flat abelian connections on a closed three-manifold.  A particular component of $\mathcal{M}_{X}$ corresponding to a bundle class $[P]\in\operatorname{Tors}H^{2}(X,\Lambda)$ is denoted as,
$\mathcal{M}_{P}\simeq H^{1}(X,\frak{t})/H^{1}(X,\Lambda)\simeq \mathbb{T}^{2g}$.  The number of components of $\mathcal{M}_{X}$ is computed for Sasakian three-manifolds in the following theorem.
\begin{theorem}\cite[Theorem 8.1]{neuray},\cite{ovz}\label{moduliprop}
Given a closed oriented Sasakian three-manifold $(X,\phi,\xi,\kappa,\operatorname{G})$ (so that $c_{1}(X)\neq 0$) then,
$$\mathcal{M}_{X}\simeq \mathbb{T}^{2g}\times \operatorname{Tors}(H^{2}(X,\Lambda))\simeq \operatorname{Hom}(\pi_{1}(X),\mathbb{T}),$$
where, $|\operatorname{Tors}H^{2}(X,\Lambda)|=|c_{1}(X)\cdot\prod_{j=1}^{M}\alpha_{j}|^{N}$.
\end{theorem}

\item  $\Omega_{P}:=\sum_{1\leq i \leq g,\\
1\leq j\leq N } d\theta_{i,j}\wedge d\bar{\theta}_{i,j}$ is the standard symplectic form on $\mathcal{M}_{P}$,

\item  $\omega_{P}:=\frac{\Omega^{gN}}{(gN)!(2\pi)^{2gN}}\in\Omega^{2gN}(\mathcal{M}_{P},\mathbb{R})$, and $\omega\in\Omega^{2gN}(\mathcal{M}_{X},\mathbb{R})$ is the symplectic form such that its restriction to the connected component $\mathcal{M}_{P}$ is $\omega_{P}$.
%$\omega_{P}:=\frac{\left(\sum_{j=1}^{g}d\theta_{j}\wedge d\bar{\theta}_{j}\right)^{gN}}{(gN)!(2\pi)^{2gN}}$.

\item  $K_{X}=\frac{1}{|c_{1}(X)\cdot \prod_{i}\alpha_{i}|^{N/2}}$,

\item  $\sqrt{T_{X}}\in\Omega^{2gN}(\mathcal{M}_{X},\mathbb{R})$ is the (square-root) of the analytic torsion (see Def. \ref{torsdef2} and Remark \ref{baserem}).  We also write $\sqrt{T_{X}}\in\Omega^{2gN}(\mathcal{M}_{P},\mathbb{R})$ when restricting $\sqrt{T_{X}}$ to a connected component $\mathcal{M}_{P}$. 

\item  The eta-invariant for the odd signature operator, $\operatorname{L}^{\operatorname{o}}$, acting on $\Omega^{1}(X,\frak{t})\oplus \Omega^{3}(X,\frak{t})$, is defined by analytic continuation,
\begin{equation}
\eta(\operatorname{L}^{\operatorname{o}}):=\lim_{s\rightarrow 0}\sum_{\lambda\in\operatorname{spec}^{*}(\operatorname{L}^{\operatorname{o}})}\operatorname{sgn}(\lambda)|\lambda|^{-s}.
\end{equation}
The eta-invariant is an analytic invariant introduced by Atiyah, Patodi and Singer \cite{aps1} defined for an elliptic and self-adjoint operator.  We note that as in \cite[Prop. 4.20]{aps1}, we may remove some spectral symmetry and the eta-invariant of $\operatorname{L}^{\operatorname{o}}$ coincides with the eta-invariant of the operator $\star d$ restricted to $\Omega^{1}(X,\frak{t})\cap \operatorname{Im}(\star d)$.

\item  $\eta_{\operatorname{grav}}(\operatorname{G})$ denotes the eta-invariant for the operator $\star d$ acting on $\Omega^{1}(X,\R)$, so that,
\begin{equation}\label{tetad}
\eta(\star d)=N\cdot\eta_{\operatorname{grav}}(\operatorname{G}),
\end{equation}
where the eta-invariant on the left hand side of \eqref{tetad} is defined on $\Omega^{1}(X,\frak{t})$ and $N=\operatorname{dim}\mathbb{T}$,

\item
\begin{equation}
\operatorname{CS}_{s}(A^{\operatorname{G}}):=\frac{1}{4\pi}\int_{X}s^{*}\operatorname{Tr}(A^{\operatorname{G}}\wedge dA^{\operatorname{G}}+\frac{2}{3} A^{\operatorname{G}}\wedge A^{\operatorname{G}}\wedge A^{\operatorname{G}}),
\end{equation}
is the gravitational Chern-Simons term, where $A^{\operatorname{G}}$ is the Levi-Civita connection and $s$ a trivializing section of twice the tangent bundle of $X$.  More explicitly, let $H=\operatorname{Spin}(6)$, $Q=TX\oplus TX$ viewed as a principal
$\operatorname{Spin}(6)$-bundle over $X$, $\operatorname{G}\in\Gamma(S^{2}(T^{*}X))$ a Riemannian metric on $X$, $\phi:Q\rightarrow \operatorname{SO}(X)$ a principal bundle morphism, and $A^{LC}\in\mathcal{A}_{SO(X)}:=\{ A\in (\Omega^{1}(\operatorname{SO}(X))\otimes\frak{so}(3))^{\operatorname{SO}(3)}\,\,|\,\,A(\xi^{\sharp})=\xi, \,\,\forall\,\xi\in\frak{so}(3)\}$ the Levi-Civita connection.  Then $A^{\operatorname{G}}:=\phi^{*}A^{LC}\in \mathcal{A}_{Q}:=\{ A\in (\Omega^{1}(Q)\otimes\frak{h})^{H}\,\,|\,\,A(\xi^{\sharp})=\xi, \,\,\forall\,\xi\in\frak{h}\}$.
\end{itemize}
An Atiyah-Patodi-Singer theorem, \cite[Prop. 4.19]{aps2}, says that the combination,
\begin{equation}\label{apsthm}
\eta_{\operatorname{grav}}(\operatorname{G})+\frac{1}{3}\frac{\operatorname{CS}(A^{\operatorname{G}})}{2\pi},
\end{equation}
is a topological invariant depending only on a two-framing of $X$.  Recall that a two-framing is a choice of a homotopy equivalence class $\Pi$ of trivializations of $TX\oplus TX$, twice the tangent bundle of $X$.  Note that $\Pi$ is represented by the trivializing section $s:X\rightarrow Q$ above.  The possible two-framings correspond to $\Z$.  The identification with $\Z$ is given by the signature defect defined by,
\begin{equation*}
\delta(X,\Pi)=\text{sign}(W)-\frac{1}{6}p_{1}(2TW,\Pi),
\end{equation*}
where $W$ is a $4$-manifold with boundary $X$ and $p_{1}(2TW,\Pi)$ is the relative Pontrjagin number associated to the framing $\Pi$ of the bundle $TX\oplus TX$.  The canonical two-framing $\Pi^{c}$ corresponds to $\delta(X,\Pi^{c})=0$.
\begin{rem}
Before we present the main quantities of interest in Definitions \ref{csdef1}, \ref{csdef2}, we note that both definitions implicitly require a choice of base $h^{0}$ for $H^{0}(X,\R)$ to be well defined.  We elaborate on this point in \S \ref{rtorssymsec}.
\end{rem}
\begin{define}\cite{mcl3}\label{csdef1}
Let $k\in\Z$ and $X$ a closed, oriented three-manifold.  The abelian Chern-Simons partition function, $Z_{\mathbb{T}}(X,k)$, is the quantity,
\begin{equation}
Z_{\mathbb{T}}(X,k)=\sum_{P\in\operatorname{Tors}H^{2}(X,\Lambda)}Z_{\mathbb{T}}(X,P,k),
\end{equation}
where,
\begin{equation}
Z_{\mathbb{T}}(X,P,k):=k^{m_X}e^{i k \operatorname{CS}_{X,P}(A_{P})}e^{\pi i N\left(\frac{\eta_{\operatorname{grav}}(\operatorname{G})}{4}+\frac{1}{12}\frac{\operatorname{CS}(A^{\operatorname{G}})}{2\pi}\right)}\int_{\mathcal{M}_{P}}\sqrt{T_{X}}.
\end{equation}
\end{define}
\begin{define}\label{csdef2}\cite{mcl3}
Let $k\in \Z$, and let $(X,\phi,\xi,\kappa,\operatorname{G})$ be a closed oriented Sasakian three-manifold.
Define the \emph{symplectic abelian Chern-Simons partition function},
\begin{equation}
\overline{Z}_{\mathbb{T}}(X,k)=\sum_{[P]\in\operatorname{Tors}H^{2}(X,\Lambda)}\overline{Z}_{\mathbb{T}}(X,P,k),
\end{equation}
where,
\begin{equation}
\overline{Z}_{\mathbb{T}}(X,P,k)=k^{m_X}e^{i k \operatorname{CS}_{X,P}(A_{P})}e^{i\pi\left(\frac{N}{4} -\frac{1}{2}\eta_{0}\right)}\int_{\mathcal{M}_{P}}K_{X}\cdot\omega_{P}.
\end{equation}
\end{define}
\noindent
The main motivation for this work is the conjectural equivalence of the rigorous topological invariants $Z_{\mathbb{T}}(X,k)$ and $\overline{Z}_{\mathbb{T}}(X,k)$.  Note that this conjecture arises simply due to the fact that the rigorous definitions of $Z_{\mathbb{T}}(X,k)$ and $\overline{Z}_{\mathbb{T}}(X,k)$ are derived from the same heuristic Chern-Simons partition function in physics.  We note that part of this conjectural equivalence is motivated by \cite{jm} which argues that $\sqrt{T_{X}}$ is proportional to a specific scalar multiple of the natural unit symplectic volume form $\omega\in\Omega^{2gN}(\mathcal{M}_{X},\R)$ by using the group structure on the moduli space $\mathcal{M}_{X}$,
\begin{equation}\label{mainconj}
\sqrt{T_{X}}=C\cdot \left(\frac{1}{\sqrt{|\operatorname{Tors}H^{2}(X,\Lambda)|}}\cdot\omega\right),
\end{equation}
where $0\neq C\in\R$.  Note that \cite{jm} works with the case where $X$ is endowed with a \emph{regular} Sasakian structure, which corresponds to a principle $\mathbb{U}(1)$ bundle over a surface \emph{without} orbifold points.  This article studies the more general case of a three-manifold $X$ that admits a Sasakian structure.  We are able to calculate the square-root of $T_{X}$ explicitly as a specific scalar multiple of a natural symplectic volume form on the moduli space $\mathcal{M}_{X}$ using a theorem of M. Rumin and N. Seshadri \cite[Theorem 5.4]{rs}.  We obtain the following,
\begin{theorem}[Main Theorem]\label{maintheorem}
Let $(X,\phi,\xi,\kappa,\operatorname{G})$ be a closed Sasakian three-manifold.  Then,
\begin{equation}
\sqrt{T_{X}}= \frac{1}{|c_{1}(X)\cdot \prod_{i}\alpha_{i}|^{N/2}}\cdot\omega.
\end{equation}
\end{theorem}
\noindent
We note that Theorem \ref{maintheorem} combined with Theorem \ref{moduliprop} leads to an explicit computation of the symplectic volume of the moduli space.  Thus, we have the following,
\begin{corollary}
Given a closed oriented Sasakian three-manifold $(X,\phi,\xi,\kappa,\operatorname{G})$, the symplectic volume of the moduli space $\mathcal{M}_{X}$ with respect to the symplectic volume form $\sqrt{T_{X}}\in\Omega^{2gN}(\mathcal{M}_{X},\R)$ is given by,
\begin{equation}
\int_{\mathcal{M}_{X}}\sqrt{T_{X}}=\sqrt{|\operatorname{Tors}H^{2}(X,\Lambda)|}=\left|c_{1}(X)\cdot\prod_{j}\alpha_{j}\right|^{N/2}.
\end{equation}
\end{corollary}
\noindent
As a consequence of Theorem \ref{maintheorem} we obtain the following verification of the above conjecture,
\begin{corollary}
Let $k\in \Z$, and let $(X,\phi,\xi,\kappa,\operatorname{G})$ be a closed oriented Sasakian three-manifold.  Then the magnitudes of $Z_{\mathbb{T}}(X,k)$ and $\overline{Z}_{\mathbb{T}}(X,k)$ agree identically,
\begin{equation}
\left|Z_{\mathbb{T}}(X,k)\right|=\left|\overline{Z}_{\mathbb{T}}(X,k)\right|,
\end{equation}
and,
\begin{equation}
\left|Z_{\mathbb{T}}(X,k)\right|=k^{m_X}\cdot \frac{\left|\sum_{[P]\in\operatorname{Tors}H^{2}(X,\Lambda)}e^{i k \operatorname{CS}_{X,P}(A_{P})}\right|}{\sqrt{|\operatorname{Tors}H^{2}(X,\Lambda)|}}.
\end{equation}
\end{corollary}

\section{Proof of the main theorem}\label{rtorssymsec}
\noindent
In this section we prove Theorem \ref{maintheorem} and compute the square root of the analytic torsion $\sqrt{T_{X}}$ as a symplectic volume form on the moduli space of flat abelian connections $\mathcal{M}_{X}$ in the case that $X$ admits a Sasakian structure.
For simplicity, we will assume $\mathbb{T}=\operatorname{U}(1)$ in this section and set $N=1$.  
\begin{rem}\label{baserem1}
The natural quantity that shows up in the symplectic abelian Chern-Simons path integral is $\omega$ multiplied by $1/|\operatorname{Vol}(I)|$, where $$I:=\{g\in\mathcal{G}_{P}|A_{P}\cdot g=A_{P}\}\simeq \operatorname{U}(1)<\mathcal{G},$$
is the isotropy subgroup of the gauge group of a given abelian connection $A_{P}\in\mathcal{A}_{P}$.  The volume of the isotropy group, $\operatorname{Vol}(I)$, requires a choice of measure on $I\simeq \operatorname{U}(1)$, which boils down to a choice of base $h^{0}$ for $H^{0}(X,\mathbb{R})$.  We recall some of the details presently.\\
\\
In our study of abelian Chern-Simons theory \cite{mcl3}, the natural invariant metric $\operatorname{H}_{\mathcal{G}}$ on the group $\mathcal{G}$ is defined in terms of the Hodge star $\star$ for the given Sasakian metric $\operatorname{G}$ on $X$,
\begin{equation}
\operatorname{H}_{\mathcal{G}}(\theta_{1},\theta_{2}):=\int_{X}\langle\theta_{1}\wedge\star \theta_{2}\rangle,
\end{equation}
where $\theta_{1},\theta_{2}\in \operatorname{Lie}\,\mathcal{G}\simeq \Omega^{0}(X,\mathbb{R})$.
Observe that $\operatorname{H}_{\mathcal{G}}$ restricted to constant functions $\theta_{1},\theta_{2}\in\mathbb{R}\subset \operatorname{Lie}\,\mathcal{G}$ is given as follows,
\begin{eqnarray*}
\operatorname{H}_{\mathcal{G}}(\theta_{1},\theta_{2})&=&\int_{X}\langle\theta_{1}\wedge\star \theta_{2}\rangle,\\
                     &=&\left(\int_{X}\star 1\right)\cdot\langle\theta_{1},\theta_{2}\rangle.
\end{eqnarray*}
We may therefore write $\sqrt{\operatorname{H}_{\mathcal{G}}}=\left(\int_{X}\star 1\right)^{1/2}$.  Now we choose the measure $\sqrt{\operatorname{H}_{\mathcal{G}}}\,d\sigma$ on $I\simeq\operatorname{U}(1)$ such that $d\sigma=d\theta/2\pi$ setting $\int_{\operatorname{U}(1)}d\sigma=1$.  Let $\mathcal{H}^{0}(X,\mathbb{R})$ denote the harmonic $0$-forms on $X$.  Note that by definition of the de Rham map $\delta_{\operatorname{dR}}^{0}:\mathcal{H}^{0}(X,\mathbb{R})\rightarrow H^{0}(X,\mathbb{R})$, this choice of measure may be viewed as a choice of base $h^{0}$ for $H^{0}(X,\mathbb{R})\simeq \operatorname{Lie}\,\operatorname{U}(1)$ such that $\delta_{\operatorname{dR}}^{0}(2\pi)=h^{0}$.  We have,
\begin{eqnarray}\nonumber
\operatorname{Vol}(I)&:=&\int_{\operatorname{U}(1)}\sqrt{\operatorname{H}_{\mathcal{G}}}d\sigma,\\\nonumber
                &=&\sqrt{\operatorname{H}_{\mathcal{G}}},\,\text{since $\int_{\operatorname{U}(1)}d\sigma=1$,}\\\label{volres}
                &=&\left[\int_{X}\star 1\right]^{1/2}.
\end{eqnarray}   
Since the Hodge star $\star$ is defined in terms of the given Sasakian metric, we have,
$$\operatorname{Vol}(I)=\left[\int_{X}\kappa\wedge d\kappa\right]^{1/2}=\left[c_{1}(X)\right]^{1/2}.$$
\end{rem}
\noindent
A proof of Theorem \ref{maintheorem} follows from \cite[Theorem 5.4]{rs}, where the analytic torsion is computed on a closed Sasakian three-manifold twisted by a unitary representation $\rho:\pi_{1}(X)\rightarrow \operatorname{U}(r)$.  Combining this with a substitution of some known special values of the Riemann-Hurwitz zeta function completes the proof. \\
\\
Let $(\operatorname{M},\operatorname{G})$ be a closed oriented Riemannian manifold of dimension $m$ and let $\rho:\pi_{1}(\operatorname{M})\rightarrow \operatorname{U}(1)$ be a representation of the fundamental group of $\operatorname{M}$.  Recall that $\rho$ corresponds to a flat principal $\operatorname{U}(1)$ bundle $P$ over $\operatorname{M}$ equipped with a flat connection $A_{\rho}\in\mathcal{A}_{P}$.  Given a representation $\chi:\operatorname{U}(1)\rightarrow \operatorname{Aut}\F$, where $\F=\R$ or $\C$, we obtain an associated line bundle $\mathcal{E}_{\chi}:=P\times_{\chi}\F$.  Let,
$$d_{A_{\rho}}^{\chi}:\Omega^{q}(\operatorname{M},\mathcal{E}_{\chi})\rightarrow \Omega^{q+1}(\operatorname{M},\mathcal{E}_{\chi}),$$
denote the covariant derivative associated to $A_{\rho}$ and let,
$$\Delta_{q}^{\chi}(\rho):=(d_{A_{\rho}}^{\chi})^{*}d_{A_{\rho}}^{\chi}+d_{A_{\rho}}^{\chi}(d_{A_{\rho}}^{\chi})^{*}:\Omega^{q}(\operatorname{M},\mathcal{E}_{\chi})\rightarrow \Omega^{q}(\operatorname{M},\mathcal{E}_{\chi}),$$
denote the corresponding Laplacian.  Define the determinant line,
\begin{equation*}
\operatorname{det}H^{\bullet}(\operatorname{M},d_{A_{\rho}}^{\chi}):=\bigotimes_{j=0}^{3}\operatorname{det}H^{j}(\operatorname{M},d_{A_{\rho}}^{\chi})^{(-1)^{j+1}},
\end{equation*}
where a superscript $-1$ denotes the dual space.  Let $|\cdot |_{L^{2}(\Omega^{\bullet}(X))}$ denote the $L^2$-metric on $\operatorname{det}H^{\bullet}(\operatorname{M},d_{A_{\rho}}^{\chi})$ induced by the identification of $H^{\bullet}(\operatorname{M},d_{A_{\rho}}^{\chi})$ with the harmonic forms $\mathcal{H}^{\bullet}(\operatorname{M},d_{A_{\rho}}^{\chi})$ via the de Rham map $\delta_{\operatorname{dR}}^{q}:\mathcal{H}^{q}(\operatorname{M},d_{A_{\rho}}^{\chi})\rightarrow H^{q}(\operatorname{M},d_{A_{\rho}}^{\chi})$.
\begin{define}\label{torsdef2}\cite{rsi}
Let $\operatorname{M}$ be a closed oriented Riemannian manifold of dimension $m$ and let $\rho:\pi_{1}(\operatorname{M})\rightarrow \operatorname{U}(1)$ be a representation of the fundamental group of $\operatorname{M}$ and let $\chi:\operatorname{U}(1)\rightarrow \operatorname{Aut}\F$ be a representation of $\operatorname{U}(1)$ (where $\F=\R$ or $\C$).  Let $\Delta_{q}^{\chi}(\rho):\Omega^{q}(\operatorname{M},\mathcal{E}_{\chi})\rightarrow \Omega^{q}(\operatorname{M},\mathcal{E}_{\chi})$ denote the Laplacian in the representation $\chi$.  Let $\zeta_{q}(s)$ be the zeta-function for $\Delta_{q}^{\chi}(\rho)$ defined for $\operatorname{Re}(s)\gg0$ by,
\begin{equation}
\zeta_{q}(s):=\frac{1}{\Gamma(s)}\int_{0}^{\infty}t^{s-1}\operatorname{tr}(e^{t\Delta_{q}}-\Pi_{q})dt,
\end{equation}
analytically continued to $\C$ and $\Pi_{q}:\Omega^{q}(\operatorname{M},\rho)\rightarrow \mathcal{H}^{q}(\operatorname{M},\rho)$ orthogonal projection.  The \emph{analytic torsion} is defined as,
\begin{equation}
T_{\operatorname{M}}=T_{\operatorname{M}}^{\chi}(\rho):=\operatorname{exp}\left(\frac{1}{2}\sum_{q=0}^{m}(-1)^{q}q\zeta_{q}'(0)\right).
\end{equation}
The \emph{Ray-Singer metric} $||\cdot ||_{RS}$ is defined as,
\begin{equation}
||\cdot ||_{RS}=T_{\operatorname{M}}\,|\cdot |_{L^{2}(\Omega^{\bullet}(X))}. 
\end{equation}
\end{define}
\noindent
Note that \cite{rs} defines and studies a new type of analytic torsion on contact manifolds called the \emph{contact analytic torsion}, denoted by $T^{C}_{X}$, and they also introduce a corresponding \emph{contact Ray-Singer metric}, denoted $||\cdot ||_{C}$.  These quantities are defined in terms of the \emph{contact complex} $(\mathcal{E},D_{H})$, originally introduced by M. Rumin \cite{rumin}, on a contact manifold $(X,\kappa)$.  Given the Reeb vector field $\xi\in\Gamma(X)$ for the contact form $\kappa\in\Omega^{1}(X,\mathbb{R})$, let $d_{H}:\Omega^{j}(X)\rightarrow \Omega^{j+1}(X)$ be defined as $d_{H}:=d-\kappa\wedge \iota_{\xi}$, and $\mathcal{L}_{\xi}$ be the Lie derivative.  Define $\Omega^{1}(H):=\{\alpha\in\Omega^{1}(X)\,|\,\iota_{\xi}\alpha=0\}$ and $\Omega^{2}(V):=\{\beta\in\Omega^{2}(X)\,|\,\beta=\kappa\wedge\alpha,\,\text{for}\, \alpha\in\Omega^{1}(X)\}$.  Given a contact metric manifold $(X, \phi, \xi, \kappa, \operatorname{G})$, and $\star$ the usual Hodge star for the metric $\operatorname{G}$, the \emph{horizontal Hodge star} is defined as $\star_{H}:=\star \circ (\kappa\wedge)$.  The \emph{contact complex} $(\mathcal{E},D_{H})$ is defined as,
\begin{equation}
C^{\infty}(X)\xrightarrow{\text{$D_{H}=d_{H}$}}\Omega^{1}(H)\xrightarrow{\text{$D_{H}=D$}}\Omega^{2}(V)\xrightarrow{\text{$D_{H}=d$}}\Omega^{3}(X),
\end{equation}
with middle operator $D_{H}=D=\kappa\wedge(\mathcal{L}_{\xi}+d_{H}\star_{H} d_{H})$.  Note that this complex may be defined using only the choice of a contact 2-plane field \cite{rs}, and we have introduced a contact metric structure in order to be more explicit.  Also note that one can twist the contact complex with a flat bundle and define the twisted contact complex, contact analytic torsion and contact Ray-Singer metric as well \cite{rs}.  Given a contact metric manifold $(X, \phi, \xi, \kappa, \operatorname{G})$, the contact analytic torsion and metric are defined using the \emph{contact Laplacian} on $(\mathcal{E},D_{H})$,
\begin{equation}
\Delta^{C}_{q}=
\begin{cases} (d_{H}^{*}d_{H}+d_{H}d_{H}^{*})^{2} & \text{if $q = 0,3$,}
\\
D^{*}D+(d_{H}d_{H}^{*})^{2} &\text{if $q=1$,}
\\
DD^{*}+(d_{H}^{*}d_{H})^{2} &\text{if $q=2$.}
\end{cases}
\end{equation}
This operator is \emph{maximally hypoelliptic and invertible in the Heisenberg symbolic calculus} \cite{rs}; a key property that allows one to make sense of the zeta function for the contact Laplacian $\zeta(\Delta_{q}^{C})(s)$.  \cite{rs} introduce the \emph{contact torsion function},
\begin{equation}
K(s):=\frac{1}{2}\sum_{q=0}^{3}(-1)^{q}w(q)\zeta(\Delta_{q}^{C})(s),
\end{equation}
where for $q=0,1,2,3$,
\begin{equation}
w(q)=\begin{cases}
q,\,\,q\leq 1,\\
q+1,\,\,q>1.
\end{cases}
\end{equation} 
Note that our definition of $K(s)$ is the negative of the one that occurs in \cite{rs}.  The \emph{contact analytic torsion} is then defined to be,
\begin{equation}
T_{X}^{C}:=\operatorname{exp}\left(\frac{1}{2}K'(0)\right).
\end{equation}       
It is shown in \cite{rs} that the analytic torsion and Ray-Singer metric agree with their contact geometric counterparts on Sasakian manifolds.  Note that our definition of $T_{X}^{C}$ is the inverse of the definition used in \cite{rs}.  
\begin{theorem}\cite[Theorem 4.2]{rs}\label{rsthm1}
Let $(X,\phi,\xi,\kappa,\operatorname{G})$ be a closed Sasakian (CR-Seifert) three-manifold, $\rho:\pi_{1}(X)\rightarrow U(N)$ a unitary representation, and $\chi_{0}:U(N)\rightarrow \operatorname{Aut}(\mathbb{C}^{N})$ the standard representation.  Let $T_{X}$ and $T_{X}^{C}$ denote the analytic torsion and the contact analytic torsion, respectively, in the standard representation; e.g. $T_{X}:=T^{\chi_{0}}_{X}$.  Then the analytic torsion $T_{X}$ and the contact analytic torsion $T^{C}_{X}$ agree,
\begin{equation}\label{rseq}
T_{X}(\rho)=T^{C}_{X}(\rho).
\end{equation}
Also, the Ray-Singer metric $||\cdot ||_{RS}$ and the contact Ray-Singer metric $||\cdot ||_{C}$ agree,
\begin{equation}\label{rseq}
||\cdot ||_{RS}=||\cdot ||_{C}.
\end{equation}
\end{theorem}
\noindent
For $a\in (0,1]$, let $\widetilde{\zeta}(s,a)=\sum_{n\in\mathbb{N}}\frac{1}{(n+a)^{s}}$ denote the Riemann-Hurwitz zeta function, and let $\widetilde{\zeta}(s):=\widetilde{\zeta}(s,1)$ denote the Riemann zeta function.  The main result that we need is given as follows. 
\begin{thm}\cite[Theorem 5.4]{rs}\label{rsthm}
Let $(X,\phi,\xi,\kappa,\operatorname{G})$ be a closed Sasakian three-manifold.  Split $\mathcal{E}_{\chi}$ into irreducibles $\mathcal{E}_{\chi}^{\theta}$.  Then the contact torsion function spectrally decomposes as,
\begin{equation}
K(s)=\sum_{\mathcal{E}_{\chi}^{\theta}}K_{\theta}(s),
\end{equation}
such that,
\begin{itemize}
\item  On $\mathcal{E}_{\chi}^{\theta}$ with $\theta\in(0,1)$, i.e. $\chi\circ\rho(h)=e^{2\pi i\theta}\neq 1$, we have,
\begin{eqnarray}
K_{\theta}(s)&=&-\operatorname{dim}(\mathcal{E}_{\chi}^{\theta})\chi(\Sigma^{*})(\widetilde{\zeta}(2s,\theta)+\widetilde{\zeta}(2s,1-\theta))\\
                                         &\,\,\,\,\,\,\,\,\,\,\,\,-&\sum_{i,j}\frac{1}{\alpha_{i}^{2s}}(\widetilde{\zeta}(2s,\theta_{i,j})+\widetilde{\zeta}(2s,1-\theta_{i,j})).
\end{eqnarray}

\item  Let $\mathcal{E}_{\chi}^{0,i}=\operatorname{ker}(1-\chi\circ\rho(c_{i}))$.  Then we have,
\begin{eqnarray*}
K_{0}(s)&=&-K(X,\rho)(2\widetilde{\zeta}(2s)+1)-2\widetilde{\zeta}(2s)\sum_{i}\operatorname{dim}(\mathcal{E}_{\chi}^{0,i})(\alpha_{i}^{-2s}-1)\\
        &\,\,\,\,\,\,\,\,\,\,-&\sum_{\{(i,j):\theta_{i,j}\neq 0\}}\frac{1}{\alpha_{i}^{2s}}(\widetilde{\zeta}(2s,\theta_{i,j})+\widetilde{\zeta}(2s,1-\theta_{i,j})).
\end{eqnarray*}
where $K(X,\rho):=2\operatorname{dim}H^{0}(X,\frak{t})-\operatorname{dim}H^{1}(X,\frak{t})$.
\end{itemize}
\end{thm}
\begin{rem}
We note that the proof of this theorem follows by application of the Riemann-Roch-Kawasaki formula \cite{rrk1}, \cite{rrk2}.
\end{rem}
\noindent
The case of interest for us is the trivial representation $\rho_{0}:\pi_{1}(X)\rightarrow \operatorname{U}(1)$.  Since this is already scalar we have,
\begin{equation}
K(s)=K_{0}(s),
\end{equation}
where, by Theorem \ref{rsthm}, we have,
\begin{equation}
K_{0}(s)=-K(X,\rho)(2\zeta(2s)+1)-2\zeta(2s)\sum_{i}(\alpha_{i}^{-2s}-1).
\end{equation}
Now we use the identification of the analytic torsion and the contact analytic torsion given in Theorem \ref{rsthm1} to write $T^{\chi}_{X}(\rho_{0})=\operatorname{exp}(K_{0}'(0)/2)$.  We compute $K_{0}'(0)$ using Theorem \ref{rsthm}.
Using the special values of the Riemann-zeta function, $\zeta(0)=-1/2$ and $\zeta'(0)=-\ln(2\pi)/2$ \cite{floyd}, and $K(X,\rho)=2\operatorname{dim}H^{0}(X,\frak{t})-\operatorname{dim}H^{1}(X,\frak{t})$ \cite[Eq. 42]{rs}, we obtain,
\begin{equation}
K_{0}'(0)/2=(2-2g)\ln(2\pi)-\sum_{i}\ln(\alpha_{i}).
\end{equation}
Thus,
\begin{equation}
T^{\chi}_{X}(\rho_{0})=\frac{(2\pi)^{2-2g}}{\prod_{i}\alpha_{i}}.
\end{equation}
It is easy to see that $T^{\operatorname{Ad}}_{X}(\rho)=T^{\chi}_{X}(\rho_{0})$ when $\rho_{0}\equiv 1$ is the trivial representation, $\chi$ is the standard representation, and $\rho:\pi_{1}(X)\rightarrow \operatorname{U}(1)$ is arbitrary.  This follows because the spectra of the corresponding Laplacians are identical.  That is, for the standard representation $\chi$, the Laplacian at the trivial representation $\rho_{0}$ is given by,
$$\Delta_{j}^{\chi}(\rho_{0}):=d^{*}d+dd^{*}:\Omega^{j}(X,\C)\rightarrow \Omega^{j}(X,\C),$$
where $d_{A_{\rho_{0}}}^{\chi}=d$ is just the ordinary de Rham derivative.  Also, for the adjoint representation,
$$\Delta_{j}^{\operatorname{Ad}}(\rho):=d^{*}d+dd^{*}:\Omega^{j}(X,\R)\rightarrow \Omega^{j}(X,\R),$$
since $d_{A_{\rho}}^{\operatorname{Ad}}=d$ for \emph{any} representation $\rho$.  Clearly, these operators have identical spectra.
By Poincar\'{e} duality $H^{3}(X,d)^{-1}$ is canonically isomorphic to $H^{0}(X,d)$, and $H^{1}(X,d)^{-1}$ is canonically isomorphic to $H^{2}(X,d)$.  Thus,
\begin{equation*}
||\cdot ||_{RS}\in |\operatorname{det}H^{0}(X,d_{A_{\rho}})|^{\otimes 2}\bigotimes|\operatorname{det}H^{1}(X,d_{A_{\rho}})^{-1}|^{\otimes 2},
\end{equation*}
and we may define the square-root of $||\cdot ||_{RS}$,
\begin{equation*}
\sqrt{||\cdot ||_{RS}}\in |\operatorname{det}H^{0}(X,d_{A_{\rho}})|\bigotimes|\operatorname{det}H^{1}(X,d_{A_{\rho}})^{-1}|.
\end{equation*}
Note that since the adjoint representation is trivial on $\mathbb{R}$, we have,
\begin{equation*}
\sqrt{||\cdot ||_{RS}}\in |\operatorname{det}H^{0}(X,\mathbb{R})|\bigotimes|\operatorname{det}H^{1}(X,\mathbb{R})^{-1}|.
\end{equation*}
\begin{rem}\label{baserem}
Observe that if $\nu^{0}$ is an orthonormal base for $\mathcal{H}^{0}(X,\R)=\R$, then it may be identified as a scalar $\nu^{0}\in\R$ such that,
\begin{eqnarray*}
1&=&||\nu^{0}||^{2},\\
 &=&\int_{X}\nu^{0}\wedge \star \nu^{0},\\
 &=&|\nu^{0}|^{2}\int_{X}\kappa\wedge d\kappa,\\
 &=&|\nu^{0}|^{2}\cdot c_{1}(X).
\end{eqnarray*}
Thus, $|\nu^{0}|=1/|c_{1}(X)|^{1/2}$.  In order to view the analytic torsion as a volume form on $\mathcal{M}_{X}$, we must choose a base $h^{0}$ for $H^{0}(X,\R)$ and evaluate $\sqrt{T_{X}}$ at $h_{0}$.  If we identify $H^{0}(X,\R)\simeq \R$ via the de Rham map $\delta_{\operatorname{dR}}^{0}$, then we make the same choice as in Remark \ref{baserem1} and choose $\delta_{\operatorname{dR}}^{0}(2\pi)=h^{0}$.  
\end{rem}
\noindent
Choosing $h^{0}\in H^{0}(X,\mathbb{R})$ as in Remark \ref{baserem} and denoting the Ray-Singer metric evaluated at $h^{0}$ by $||\cdot ||_{RS}|_{h^{0}}$, we define,
\begin{equation}
\sqrt{T_{X}}:=\sqrt{||\cdot ||_{RS}|_{h^{0}}}\in |\operatorname{det}H^{1}(X,\mathbb{R})^{-1}|.
\end{equation}
We therefore have,
\begin{equation}
\sqrt{T_{X}}=\frac{(2\pi)^{-Ng}}{|c_{1}(X)\cdot \prod_{i}\alpha_{i}|^{N/2}}\left|\bigwedge \delta_{\operatorname{dR}}^{1}(\nu^{1})\right|^{*},
\end{equation}
where $\left|\bigwedge \delta_{\operatorname{dR}}^{1}(\nu^{1})\right|^{*}:\bigwedge^{\operatorname{max}}H^{1}(X,\frak{t})\rightarrow \R^{+}$ is the volume form associated to the basis given by $\delta_{\operatorname{dR}}^{1}(\nu^{1})$.  
Writing the above results concisely, if $(X,\phi,\xi,\kappa,\operatorname{G})$ is a closed Sasakian three-manifold, then,
\begin{equation}
\sqrt{T_{X}}= \frac{1}{|c_{1}(X)\cdot \prod_{i}\alpha_{i}|^{N/2}}\cdot\omega,
\end{equation}
where,
\begin{equation}\label{omegaeqn}
\omega:=\frac{\Omega^{gN}}{(gN)!(2\pi)^{2gN}},
\end{equation}
and,
\begin{equation}
\Omega:=\sum_{1\leq i \leq gN } d\theta_{i}\wedge d\bar{\theta}_{i}.
\end{equation}
Note that the generalization to the case of an arbitrary torus $\mathbb{T}$ is straightforward.  We also point out that the extra factor of $(2\pi)^{gN}$ that occurs in Eq. \eqref{omegaeqn} is due to the corresponding factor of $\sqrt{2\pi}$ in the norm of each orthonormal basis element for the first cohomology.

\section*{Acknowledgement}
\begin{itemize}
\item  Supported in part by the Danish National Research Foundation grant DNRF95 (Centre for Quantum Geometry of Moduli
Spaces – QGM).

\item  Supported in part by the European Science Foundation Network ‘Interactions of Low-Dimensional Topology and Geometry
with Mathematical Physics’.  

\item  Supported in part by a Research Instructorship at Northeastern University, Department of Mathematics, 360 Huntington Ave, Boston, MA, USA, 02115.
\end{itemize}
The author would like to thank Chris Beasley, Lisa Jeffrey and Andrew Swann for several helpful comments related to this work.    
%\biboptions{authoryear}
\bibliographystyle{elsarticle-num}
%\bibliography{finalthesis}



%\begin{thebibliography}{00}



% \end{thebibliography}

\end{document}